\documentclass[10pt]{amsart}
\usepackage{amsthm, amssymb, amsmath, graphicx, subfig}

\newtheorem{define}{Definition}
\newtheorem{note}{Note}

\def\sign{\mathrm{sign}}

\begin{document}

\title{The Concordance Genus of 11--crossing Knots}
\author{M. Kate Kearney}

 \begin{abstract} The concordance genus of a knot is the least genus of any knot in its concordance class.  It is bounded above by the genus of the knot, and bounded below by the slice genus, two well-studied invariants.  In this paper we consider the concordance genus of 11--crossing prime knots.  This analysis resolves the concordance genus of 533 of the 552 prime 11--crossing knots.  The appendix to the paper gives concordance diagrams for 59 knots found to be concordant to knots of lower genus, including null-concordances for the 30 11--crossing knots known to be slice.
 
 \end{abstract}

\maketitle

\section{Introduction}\label{intro}

The knot concordance group, defined by Fox and Milnor in 1966 \cite{FM:CobordismOfKnots}, has been a central focus of the field of knot theory for many years.  The foremost objectives in the field are algebraically, to understand the structure of the concordance group, and topologically, to understand the relationship between knots in a given concordance class.  It is invaluable to this process to consider examples constructed to satisfy certain properties.  Our focus is in furthering the development of examples by looking at low crossing number prime knots from the perspective of the concordance.  In particular, to consider the relationship  of the concordance genus to several other invariants, and to use it as a computational tool to find concordances between low crossing number knots.

The concordance genus is the least genus of a surface in $S^3$ whose boundary is a knot concordant to the knot.  Recall that two knots are concordant if they co-bound an annulus embedded in $S^3 \times I$.  Although a variety of invariants provide bounds on the concordance genus, the actual determination can be difficult even for low crossing number knots.  It is immediately evident that $g_4(K) \leq g_c(K) \leq g_3(K)$ (where $g_4$ denotes the smooth four-genus, or slice genus, $g_c$ denotes the smooth concordance genus, and $g_3$ denotes the genus of the knot).  Casson gave a first example of a knot for which $g_4(K) \neq g_c(K)$, and Nakanishi showed that this gap could be made arbitrarily large \cite{N:Unknotting}.  For each positive integer $n$, Gilmer gave examples of (algebraically slice) knots with $g_4(K) = g_3(K) = n$ (and hence the concordance genus is also $n$) \cite{G:SliceGenus}.  The Alexander polynomial and signature of knots provide further bounds on the concordance genus as will be discussed in this paper.  Slice genus bounds in particular have been studied in depth by a number of topologists.  The bounds used in this paper follow from work of Murasugi, Levine, and Tristram \cite{Le:Invariants}, \cite{M:Signature}, \cite{T:Signature}.  There are other known bounds on the slice genus, such as the $\tau$ invariant of Heegaard-Floer theory and the $s$ invariant from Khovanov homology, which may help in future progress on calculations of the concordance genus as well.

A complete census of the concordance genus of low crossing number knots was begun by Livingston in \cite{L:ConcordanceGenusI} and extended to include all primes knots of $10$ or fewer crossings in \cite{L:ConcordanceGenusII}.  In this paper we consider the concordance genus of 11--crossing knots.  For 474 of the 552 prime knots of 11 crossings we observe that known invariants obstruct the knot from being concordant to a knot of lower genus.  For another 59 knots we show an explicit concordance, illustrated in the appendix.  This extends the explicit construction of slice disks begun by Kawauchi in \cite{Kawauchi}.  This paper determines the concordance genus for all but 19 of the 552 prime knots of eleven crossings.  We will be working in the smooth category throughout this paper, although many of the tools (including the signature and Alexander polynomial bounds) also apply to the topological locally flat case.

In Section 2 of this paper we give a detailed summary of our findings.  Section 3 gives background material on the obstructions used and details their relationship to the concordance genus.  In Section 4 we illustrate with several examples the process of combining the obstructions with geometric constructions to determine the concordance genus.

\section{Summary of Results}\label{summary}

The calculations of the concordance genus for prime 11--crossing knots can be broken into two general strategies, finding a concordance to a knot of lower genus, and finding obstructions that prevent a concordance to a knot of lower genus.  For the majority of prime 11--crossing knots, bounds from classical invariants are sufficient to determine the concordance genus.  In these cases we confirm that the concordance genus is equal to the three-genus.  Slice knots are concordant to the unknot, and therefore have concordance genus 0.  For many of the remaining knots, the concordance genus can be determined by finding a concordance to a knot of lower genus.  There are several knots for which the concordance genus has not yet been determined.

Two of the most useful obstructions are given by the Alexander polynomial and the signature of the knot.  We begin by recalling that the three-genus is an upper bound for the concordance genus, so with the help of a lower bound we can determine the concordance genus in some cases.  If the Alexander polynomial is irreducible, by the Fox-Milnor theorem, the degree of the Alexander polynomial is a lower bound for twice the concordance genus, as discussed in Section \ref{polynomial_and_signature}.  With this information alone, the concordance genus is determined for 384 of the 552 11--crossing knots.  For another 84 knots, the Alexander polynomial can be factored, but each symmetric factor only appears once, so the degree of the Alexander polynomial bounds twice the concordance genus and shows that the concordance genus is equal to the three-genus.  For six knots, the signature of the knot is also necessary to determine that the concordance genus is equal to the three-genus.  Examples of the use of these obstructions are given in Section \ref{obstruction_examples}.  The values of these invariants and the concordance genus of these knots can be found on the knotinfo website \cite{knotinfo}.

Of the 78 remaining knots, 30 knots are known to be slice, and thus have concordance genus 0. The slice 11--crossing knots are:

\begin{align*}
	&11a_{28}, 11a_{35}, 11a_{36}, 11a_{58}, 11a_{87}, 11a_{96}, 11a_{103}, 11a_{115}, 11a_{164}, \\
	&11a_{165}, 11a_{169}, 11a_{201}, 11a_{316}, 11a_{326}, 11n_{4}, 11n_{21}, 11n_{37}, 11n_{39}, \\ 
	&11n_{42}, 11n_{49}, 11n_{50}, 11n_{67}, 11n_{73}, 11n_{74}, 11n_{83}, 11n_{97}, 11n_{116},  \\
				&11n_{132}, 11n_{139}, 11n_{172}
\end{align*}	

The sliceness of 11--crossing knots was determined by a variety of people, as referenced in detail in \cite{knotinfo}.  Diagrams for the null-concordances are given in the appendix.  
	
 There are 29 knots which are not slice, but are concordant to a knot of lower genus, as listed below, in Table \ref{tab:concordances}.

\begin{table}[h]
	\centering
		\begin{tabular} {| c || c |}
			\hline
			$3_1$ & $11a_{196}, 11a_{216}, 11a_{283}, 11a_{286}, 11n_{106}, 11n_{122}$  \\
			\hline
			$4_1$ & $11a_{5}, 11a_{104}, 11a_{112}, 11a_{168}, 11n_{85}, 11n_{100}$ \\
			\hline
			$5_1$ & $11n_{69}, 11n_{76}, 11n_{78}$ \\
			\hline
			$5_2$ & $11n_{68}, 11n_{71}, 11n_{75}$ \\
			\hline
			$6_2$ & $11a_{57}, 11a_{102}, 11a_{139}, 11a_{199}, 11a_{231}$ \\
			\hline
			$6_3$ & $11a_{38}, 11a_{44}, 11a_{47}, 11a_{187}$ \\
			\hline
			$3_1\#4_1$ & $11a_{132}, 11a_{157}$ \\
			\hline
		\end{tabular}
		\caption{Concordances}
		\label{tab:concordances}
\end{table}
The appendix gives diagrams of these concordances.  The examples in Section \ref{concordance_examples} explain the notation of the diagrams in the appendix and details of these calculations.  

\begin{table}[h]
	\centering
		\begin{tabular} {| c || c | c | c | c | c | c |}
			\hline
			Knot & $g_c$ bounds & $g_3$ & $g_4$ & $\sigma$ & $\Delta_K$ bound & possible concordance \\
			\hline
			$11a_6$ & [2,3] & 3 & [1,2] & 2 & 2 & $3_1 \# 4_1$ \\
			\hline
			$11a_8$ & [2,3] & 3 & 1 & 0 & 2 & $6_3$ \\
			\hline 
			$11a_{67}$ & [1,3] & 3 & [1,2] & 0 & 1 & $4_1$ \\
			\hline
			$11a_{72}$ & [2,4] & 4 & [1,2] & 0 & 2 & ? \\
			\hline
			$11a_{108}$ & [2,4] & 4 & [1,2] & 2 & 2 & $6_2$ \\
			\hline
			$11a_{109}$ & [2,4] & 4 & [1,2] & 0 & 2 & $6_2$ \\
			\hline
			$11a_{135}$ & [2,3] & 3 & [1,2] & 0 & 2 & ? \\
			\hline
			$11a_{181}$ & [2,3] & 3 & [1,2] & -2 & 2 & $6_2$ \\
			\hline
			$11a_{249}$ & [2,3] & 3 & [1,2] & 0 & 2 & $6_3$ \\
			\hline
			$11a_{264}$ & [2,4] & 4 & 1 & -2 & 2 & $3_1 \# 4_1$ \\
			\hline
			$11a_{297}$ & [1,3] & 3 & [1,2] & 2 & 1 & $5_2$ \\
			\hline
			$11a_{305}$ & [2,4] & 4 & [1,2] & 2 & 2 & $3_1 \# 4_1$ \\
			\hline
			$11a_{332}$ & [2,4] & 4 & [1,2] & 0 & 2 & $7_7$ \\
			\hline
			$11a_{352}$ & [2,3] & 3 & [1,2] & -2 & 2 & $3_1 \# 4_1$ \\
			\hline
			$11n_{34}$ & [0,3] & 3 & [0,1] & 0 & 0 & slice? \\
			\hline
			$11n_{45}$ & [1,3] & 3 & 1 & 0 & 0 & ? \\
			\hline
			$11n_{66}$ & [1,3] & 3 & [1,2] & -2 & 1 & $3_1$ \\
			\hline
			$11n_{145}$ & [1,3] & 3 & 1 & 0 & 0 & ? \\
			\hline
			$11n_{152}$ & [2,3] & 3 & 1 & -2 & 2 & $8_6$ \\
			\hline
		\end{tabular}
		\caption{Unknown values}
		\label{table:unknown}
\end{table}

For 19 knots the value of the concordance genus is still unknown.  Some bounds are known for these knots, although they don't determine the concordance genus.  The known information is summarized in Table \ref{table:unknown}.  We focus on the bounds that have been useful in calculating the concordance genus of other knots.  The values for $g_c$ bounds are given by the maximal lower bound based on the bounds from the four-genus, the signature, and the Alexander polynomial, and the upper bound given by the three-genus.  The Alexander polynomial and the Fox-Milnor theorem limit the list of possible knots concordant to a given knot, as described in the example in Section \ref{concordance_examples}.  For each knot in the Table \ref{table:unknown}, the knot of smallest crossing number which is potentially concordant to the given knot is listed.  In each case, an actual concordance has not yet been found.  

For the majority of the 19 knots above, the four-genus of the knot is also unknown.  In some cases the determination of the four-genus could lead to a better lower bound for the concordance genus, although in most cases the bound from the Alexander polynomial is stronger.  There are no knots for which the determination of the four-genus would lead directly to the determination of the concordance genus.  A particularly interesting knot in this list is $11n_{34}$, which is topologically slice (since it has Alexander polynomial 1), but it is not known whether it is smoothly slice.  In fact, it is a mutant of $11n_{42}$ which is smoothly slice, but as shown by Kearton \cite{K:Mutation} mutation does not preserve concordance class, so $11n_{34}$ is not necessarily smoothly slice.

\section{Background}\label{background}

In this section we will give several useful definitions and theorems.  We present the basic tools used for finding concordances or obstructing concordance to a knot of lower genus.

\subsection{Alexander Polynomial and Signatures}\label{polynomial_and_signature}

The Alexander Polynomial and Levine-Tristram signatures can both be defined in terms of the Seifert form.  A detailed account of both definitions can be found in \cite{Rolfsen}.  For a knot, $K$, in $S^3$, and a surface $F$ in $S^3$ with $\partial F = K$, we let $V$ be the matrix of the associated linking form.  That is, for $x_1, x_2, ... , x_n$ a symplectic basis for $H_1(F)$, and for $x_i^+$, the positive push-off of $x_i$, we let $V$ be the matrix given by $v_{i,j} = lk(x_i, x_j^+)$.  The Seifert matrix, $V$, itself is not a knot invariant, although up to a particular stabilization, the isomorphism class of the bilinear form given by $V$ is in fact a knot invariant (for details on this fact, refer to \cite{Rolfsen}).  There are also a number of algebraic invariants that can be determined from $V$.  Of particular interest are the Alexander polynomial and the Levine-Tristram signatures.

\begin{define}
 The Alexander Polynomial is given by $\Delta_K(t) := \det(V-tV^T)$.  It is well-defined up to multiplication by $\pm t^k$.  This equivalence will be denoted by $\doteq$.
\end{define}
\begin{define}
 For each unit complex number $\omega$ which is not a root of the Alexander Polynomial, there is a Levine-Tristram signature defined to be $\sigma_\omega(K) = \sign((1-\omega)V + (1- \bar{\omega})V^T)$, where $\sign$ denotes the algebraic signature of a Hermitian matrix.  We extend this definition to roots of the Alexander by taking the average of the limits, $$\sigma_\omega(K) := \frac{1}{2}(\lim_{\omega^+ \rightarrow \omega} \sigma_{\omega^+}(K) + \lim_{\omega^- \rightarrow \omega} \sigma_{\omega^-}(K)),$$ where $\omega^+$ and $\omega^-$ denote unit complex values approaching $\omega$ from opposite directions.  In this way we ensure that $\sigma_\omega$ is a well-defined concordance invariant.  \end{define}

\begin{define}
In the case of $\omega = -1$, the signature $\sigma_\omega (K)$ is referred to as the signature or the Murasugi signature, and will be denoted by $\sigma(K)$.

\end{define}

The Alexander polynomial and signatures are both easily computable and provide bounds on other knot invariants.  In particular, we will use the following facts:

\begin{itemize}
	\item $g_c(K) \leq g_3(K)$
	\item $\frac{1}{2}|\sigma_{\omega}(K)| \leq g_4(K) \leq g_c(K)$ \cite{Le:Invariants}, \cite{M:Signature}, \cite{T:Signature}
	\item $\frac{1}{2}\deg(\Delta_K(t)) \leq g_3(K)$
	\item For slice $K$, $\Delta_K(t) \doteq f(t)f(t^{-1})$ for some polynomial $f(t)$ (Fox-Milnor) \cite{FM:CobordismOfKnots}
\end{itemize}

As a consequence of the Fox-Milnor theorem, in certain cases, half the degree of the Alexander polynomial also bounds $g_c(K)$, as discussed in Section \ref{calculations}.

The details of the applications of these results will be further illustrated in the examples given in Section \ref{examples}.

\subsection{Slice Knots and Concordance}\label{slice_definition}

A knot is smoothly slice (or just {\sl slice} when the category is understood) if it is the boundary of a disk smoothly embedded in the four ball, such that $\partial B^2 \subset \partial B^4 = S^3$.  Two knots $K$ and $J$ are concordant if there is a smoothly embedding of $S^1 \times I$ in $S^3 \times I$ with boundary $K \cup -J$ embedded so that $K \subset S^3 \times \left\{0\right\}$ and $-J \subset S^3 \times \left\{1\right\}$.  Equivalently, $K$ and $J$ are concordant if $K \# - J$ is slice.  Both definitions give useful perspectives, so we will use them interchangeably.  It follows immediately from the first definition that concordance is an equivalence relation.  This confirms that the concordance genus is well-defined.  

\begin{note}
For the purposes of calculating the concordance genus, we do not need to distinguish between $K$ and $-K$.  Since $g_3(K)=g_3(-K)$, the concordance genus also does not distinguish $K$ from $-K$.  For this reason the distinction between a knot and its mirror is omitted in the enumeration of concordances.  The correct choice can easily be determined by the reader by performing the designated band moves.
\end{note}

One possibility for calculating the concordance genus of a knot is to find a knot concordant to the knot in question, for which the concordance genus is known.  As a special case of this, the concordance genus of any slice knot is zero.  If $K$ is slice, then $K \# -U$ is also slice (where $U$ denotes the unknot), so any slice knot is concordant to the unknot and thus $g_c(K) = g_c(U) = 0$.  There are $30$ prime 11--crossing knots which known to be slice.  The only 11--crossing knot for which sliceness is unknown is $11a_{34}$.

\subsection{Finding Concordances and Obstructions}\label{calculations}

To calculate the concordance genus, we must either find a concordance or find an obstruction to the existence of a concordance to a knot of lower genus.  In most cases for prime 11--crossing knots, a combination of obstructions are enough to determine the concordance genus.  The four-genus and half the signature give lower bounds for the concordance genus, and the three-genus gives an upper bound.  

\begin{note}
The signature actually provides a lower bound for the four-genus, and consequently for the concordance genus.  So in the case where the four-genus is known, it always gives at least as sharp a bound as the signature.  In fact, the four-genus is known for many of the prime 11--crossing knots, and for those with unknown four-genus, the bound from the Alexander polynomial bound is sharper than the signature bound.  There are some knots for which a stronger bound can be attained using the signature function and Alexander polynomial in conjunction by noting that the signature function has jumps only at roots of the Alexander polynomial \cite{Lickorish}.  An example of this is given in Section \ref{obstruction_examples}.
\end{note}

The degree of the Alexander polynomial gives a lower bound for $2g_3(K)$.  Using this along with the Fox-Milnor theorem, we can find a lower bound for the concordance genus.  We write $\Delta_K(t) \doteq g(t)f(t)f(t^{-1})$ for some polynomials $g(t)$ and $f(t)$, where $f(t)$ has maximal degree for decompositions of this form.  Then for any $K'$ concordant to $K$, $K \# -K'$ is slice, so $\Delta_{K \# -K'}(t) \doteq \Delta_K(t) \cdot \Delta_{K'}(t) \doteq h(t)h(t^{-1})$ for some $h(t)$.  Then $g(t^{-1})$ must be a factor of $\Delta_{K'}(t)$, and thus the degree of $g(t)$ (which is equal to the degree of $g(t^{-1})$) is a lower bound for $2g_c(K)$.  This proved to be a particularly valuable obstruction in the case of prime 11--crossing knots.  There were 384 knots with irreducible polynomials, and 84 which factored but had no factors of the form $f(t)f(t^{-1})$, each of these cases satisfying $\deg(\Delta_K(t)) = 2g_3(K)$.  We conclude then that $g_c = g_3$.

There were six knots for which the concordance genus could not be determined by the Alexander polynomial obstruction, but the four-genus and signatures provide a strong enough lower bound to determine the concordance genus.

For the remainder of the knots, a concordance to a knot of lower genus was found.  This task was completed by performing band moves on a knot diagram in hopes of relating the diagram to a simpler diagram.  Consider $K \times I$, then one oriented band move corresponds to adding a one handle to this surface.  If this results in an unlinked diagram of $K' \cup U$, we attach a disk along the unlinked boundary component to find a surface $F \cong S^1 \times I$ with boundary components $K$ and $K'$.  That is, $K$ is concordant to $K'$.  Similarly, we can perform several band moves, increasing the number of components with each band move, to create a surface with boundary $K' \cup U \cup U \cdots \cup U$.  If this surface has several unlinked unknotted boundary components, we can attach disks along each unknot to get a concordance between $K$ and $K'$.

It is valuable in looking for concordances to take advantage of the Fox-Milnor theorem to determine possible knots to find a concordance.  For some simple knots, $K'$,  it is useful to consider $K \# -K'$ and look for a slicing disk.  In the particular cases of the trefoil and the figure-eight knot, this can be accomplished by a clasp change in the diagram of $K$, similar to that used by Tamulis \cite{Ta:ConcordanceOrder}.

\begin{figure}[h]
		\centering
		\includegraphics[width=.7\textwidth]{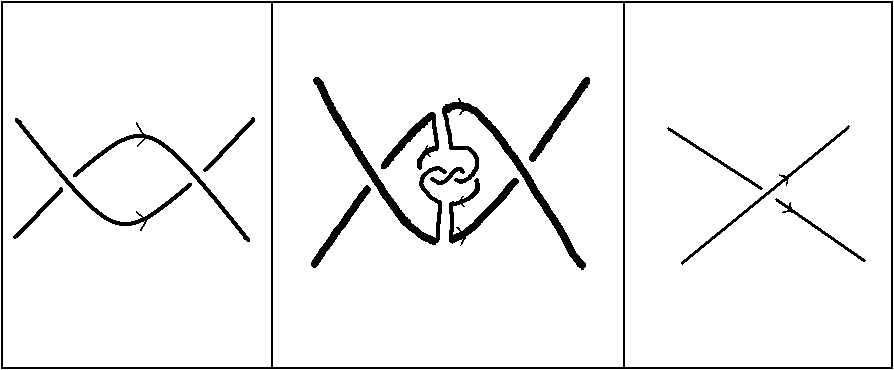}
		\caption{Clasp move for concordance to the trefoil}
	\label{fig:3_1clasp}
\end{figure}

The clasp move for the trefoil consists of taking the connected sum $K \# 3_1$, and performing one band move.  This is equivalent to the clasp change indicated in Figure \ref{fig:3_1clasp}.  This will result in two components.  If the result is the unlink of two components then $K$ is concordant to $3_1$.

\begin{figure}[h]
		\centering
		\includegraphics[width=.7\textwidth]{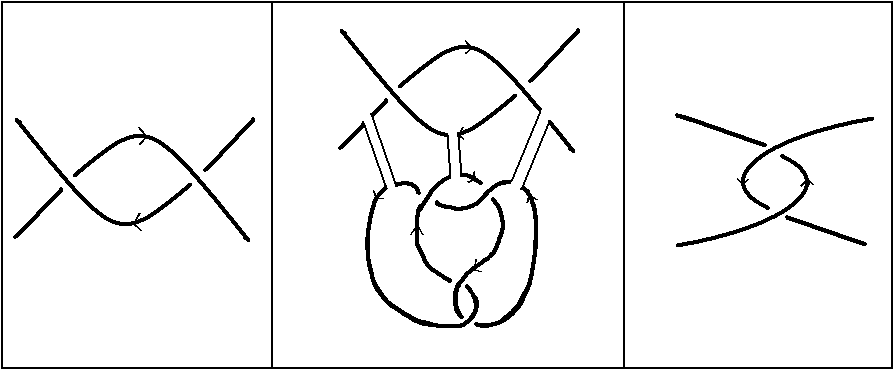}
		\caption{Clasp move for concordance to the figure eight knot}
	\label{fig:4_1clasp}
\end{figure}

The clasp move for the figure eight knot consists of taking the connected sum $K \# 4_1$, performing two band moves, and capping an unlinked component. This is equivalent to the clasp change indicated in Figure \ref{fig:4_1clasp}.  The band moves create a link with three components, one of which will be unlinked regardless of $K$.  After reducing the diagram by capping the unlinked component, the result has two components.  If this is the unlink of two components then $K$ is concordant to $4_1$.

There were 29 knots found to be concordant to knots of lower genus, as listed in Section 2 of this paper.  Diagrams for these concordances are given in the Appendix.

\section{Detailed Examples}\label{examples}

\subsection{Obstructions}\label{obstruction_examples}

For many of the prime 11--crossing knots the degree of the Alexander polynomial gives a bound equal to the three-genus.  The knot $11a_{1}$ has three-genus 3 and Alexander polynomial $\Delta_{11a_{1}}(t) \doteq 2-12t+30t^2-39t^3+30t^4-12t^5+2t^6$, which is irreducible over $\mathbb{Z}[t,t^{-1}]$.  Thus 

$$3 = \frac{1}{2}\deg(\Delta_{11a_{1}}(t)) \leq g_c(11a_{1}) \leq g_3(11a_{1}) = 3.$$

Some knots have reducible Alexander polynomials, but the degree of the Alexander polynomial still obstructs a concordance to a knot of lower genus.  The knot $11a_{51}$ has three-genus 3 and Alexander polynomial $\Delta_{11a_{51}}(t) \doteq 1-9t+28t^2-39t^3+28t^4-9t^5+t^6 \doteq (1-3t+t^2)(1-3t+5t^2-3t^3+t^4)$.  In particular, notice that there are no factors of $\Delta_{11a_{51}}(t)$ of the form $f(t)f(t^{-1})$, so for any knot, $K$, concordant to $11a_{51}$ we know that $\Delta_{11a_{51}}(t)$ is a factor of $\Delta_K(t)$.  Then for any such $K$, $\deg(\Delta_K(t)) \geq 6$.  So 

$$3 = \frac{1}{2}\deg(\Delta_{11a_{51}}(t)) \leq g_c(11a_{51}) \leq g_3(11a_{51}) = 3.$$

There are several 11--crossing knots for which the Alexander polynomial bound is not enough to determine the concordance genus, but the signature is an effective bound.  For $11a_{43}$, the Alexander polynomial is 
\begin{align*}
\Delta_{11a_{43}}(t) &\doteq 4-14t+30t^2-37t^3+30t^4-15t^5+4t^6 \\
&\doteq (1-t+t^2)^2(4-7t+4t^2),
\end{align*}
 so we can only conclude that $4-7t+4t^2$ is a factor of $\Delta_K(t)$ for any $K$ concordant to $11a_{43}$.  We can't immediately determine the concordance genus from the Alexander polynomial, but the signature is $-6$ so we have 
 $$3 = \frac{1}{2}|\sigma(11a_{43})| \leq g_4(11a_{43}) \leq g_c(11a_{43}) \leq g_3(11a_{43}) = 3.$$  
 
 On the other hand, consider $11n_{81}$: 
\begin{align*}
\Delta_{11n_{81}}(t) &\doteq 	1-3t+4t^2-4t^3+3t^4-4t^5+4t^6-3t^7+t^8 \\
&\doteq (1-t+t^2)^2(1-t+t^2-t^3+t^4).
\end{align*}
 Since the factor $(1-t+t^2)^2$ of the Alexander polynomial is already of the form $f(t)f(t^{-1})$ we may immediately only conclude that the factor $(1-t+t^2-t^3+t^4)$ must divide the Alexander polynomial of any knot concordant to $11n_{81}$ (so $g_c(11n_{81}) \geq 2$).  Since $g_3(11n_{81}) = 4$, we cannot yet determine the concordance genus.  We continue by looking at the signature.  The signature of $11n_{81}$ is $-6$, which still does not provide a sufficient bound.  However, we observe that at the roots of $(1-t+t^2)^2$ the signature function jumps by $4$.  Since the signature function is a concordance invariant and can only have jumps at roots of the Alexander polynomial, we conclude that any knot concordant to $11n_{81}$ must have $(1-t+t^2)$ as a factor of the Alexander polynomial.  Further, to satisfy the Fox-Milnor theorem, we in fact have that $(1-t+t^2)^2$ must be a factor of the Alexander polynomial of any knot concordant to $11n_{81}$.  Then for any knot, $K'$, concordant to $11n_{81}$, $4 \leq \frac{1}{2} \deg(\Delta_{K'}(t)) \leq g_3(K')$.  So we may conclude $$4 \leq g_c(11n_{81}) \leq g_3(11n_{81}) = 4.$$
 
 For each of the knots $$11a_{43}, 11a_{263}, 11n_{72}, 11n_{77}, 11n_{81}, 11n_{164}$$ the signatures and four-genus contributed to obstructing a concordance.

\subsection{Concordances}\label{concordance_examples}

\begin{figure}[h]
		\centering
		\subfloat[$11a_{196}$]{\includegraphics[width=.5\textwidth]{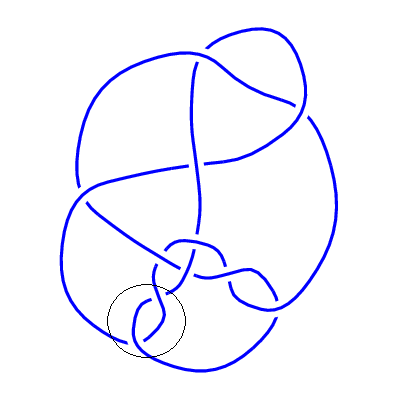}}
		\subfloat[After a clasp move, $U \# U$]{\includegraphics[width=.5\textwidth]{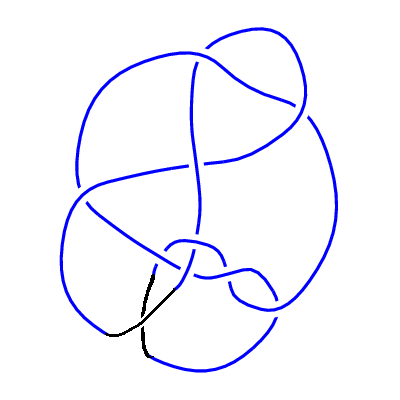}}
		\caption{$11a_{196}$ is concordant to $3_1$}
	\label{fig:ex3_1}
\end{figure}

In many examples, the bounds do not give enough information to determine the concordance genus of the knot.  We then must consider concordances.  Fortunately the algebraic invariants suggest a specific possible concordance.  For example, the Alexander polynomial of $11a_{196}$ is 
	\begin{align*}
	\Delta_{11a_{196}}(t) \doteq& 1-6t+17t^2-31t^3+37t^4-31t^5+17t^6-6t^7+t^8 \\
				\doteq& (1-t+t^2)(-1+2t-3t^2+t^3)(-1+3t-2t^2+t^3),
	\end{align*}
 so by the Fox-Milnor theorem \cite{FM:CobordismOfKnots}, if $11a_{196}$ is concordant to another knot, $K$, then $1-t+t^2$ must factor $\Delta_K(t)$.  We know that $\Delta_{3_1}(t) \doteq 1-t+t^2$, so we consider the possibility that $11a_{196}$ is concordant to $3_1$.  The concordance shown in Figure \ref{fig:ex3_1} confirms that this is true.  Note that the circle in this diagram indicates using the clasp move for $3_1$ described in Figure 1 in Section 3.3.  After changing the clasp, the diagram is two unlinked circles.  We cap off each circle with a disk to get a slicing disk for $11a_{196} \# 3_1$.

\begin{figure}[h]
		\centering
		\subfloat[$11a_{57} \# 6_2$]{\includegraphics[width=.5\textwidth]{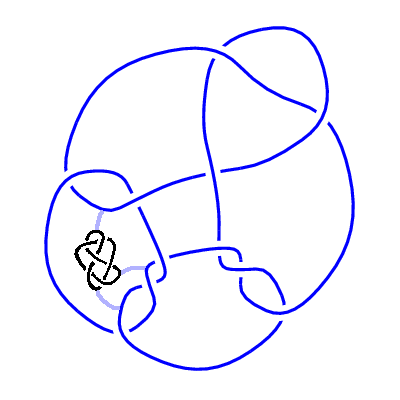}}
		\subfloat[After band moves, $U \cup U \cup U$]{\includegraphics[width=.5\textwidth]{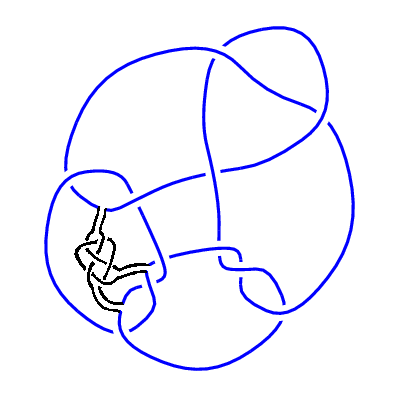}}
		\caption{$11a_{57}$ is concordant to $6_2$}
	\label{fig:ex6_2}
\end{figure}

\begin{figure}[h]
		\centering
		\subfloat[$11n_{69}$]{\includegraphics[width=.5\textwidth]{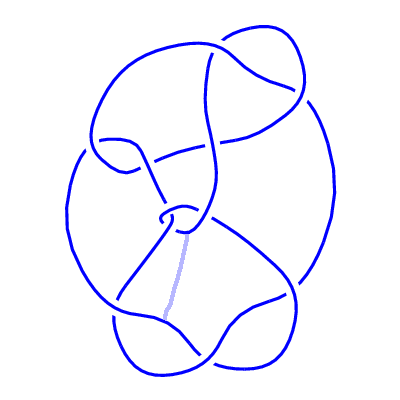}}
		\subfloat[After a band move, $5_1 \cup U$]{\includegraphics[width=.5\textwidth]{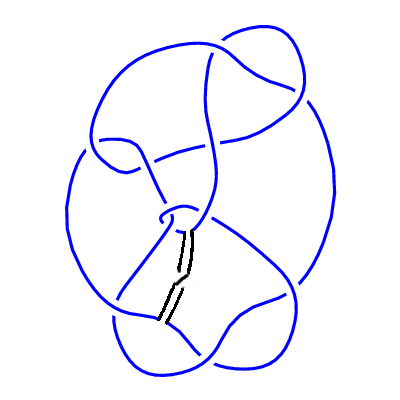}}
		\caption{$11n_{69}$ is concordant to $5_1$}
	\label{fig:ex5_1}
\end{figure}

Depending on the knot, a concordance may be easier to find in either of two ways: by showing that $K\#K'$ is slice (as in Figure \ref{fig:ex6_2}), or by directly finding a cylinder with boundary $K \cup K'$ (as in Figure \ref{fig:ex5_1}).  In either case the goal is similar.  Consider $K \times I$ as a surface with boundary $K \cup -K$.  Then add one or more bands ($D^1 \times I$) to this surface by an orientable attaching map and observe the new boundary, as in Figure \ref{fig:ex5_1}.  If the modified boundary component of this surface is a union of unlinked unknots, we can cap each off with a disk to build a slicing disk for $K$.  If the modified boundary is a union $K'$ and several unlinked unknots, we similarly cap off each disk and we have a concordance to $K'$.  Note that for each band added we should gain one extra unlinked component.  If we decrease the number of components, the surface we have built is no longer diffeomorphic to $S^1 \times I$, and consequently does not give a concordance.

\appendix
\section{Concordance Diagrams}\label{diagrams}

Included in this appendix are diagrams of all the new concordances found through this project, which are all of the known concordances for 11--crossing prime knots.  They are grouped by concordance classes.  The diagrams include the knot diagram of the 11--crossing knot, along with notation indicating how to find a concordance to the given knot.  A circle indicates use of the clasp switch for concordances to $3_1$ or $4_1$ as described in Section 3.3.  A grey arc indicates that an oriented band move is used, also described in Section 3.3.  In these diagrams the choice of orientation and twist of the band is the obvious choice which matches the orientations of the arcs it connects.  The examples in Section 4.2 give further detail about the development and notation of these diagrams for a few knots.

\begin{figure}[h]
		\centering
		\subfloat[$11a_{196}$]{\includegraphics[width=5cm]{11a_196.png}}
		\subfloat[$11a_{216}$]{\includegraphics[width=5cm]{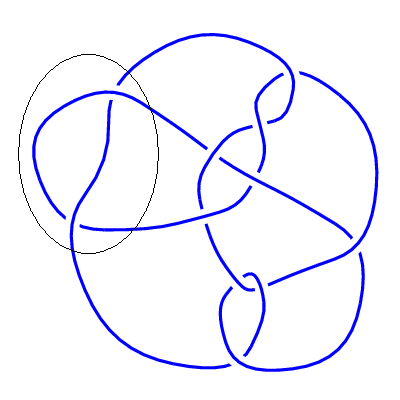}}
		
		\subfloat[$11a_{283}$]{\includegraphics[width=5cm]{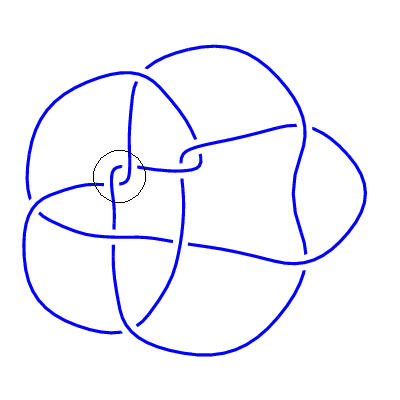}}
		\subfloat[$11a_{286}$]{\includegraphics[width=5cm]{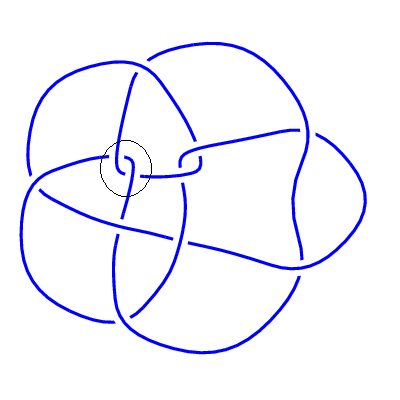}}
		
		\subfloat[$11n_{106}$]{\includegraphics[width=5cm]{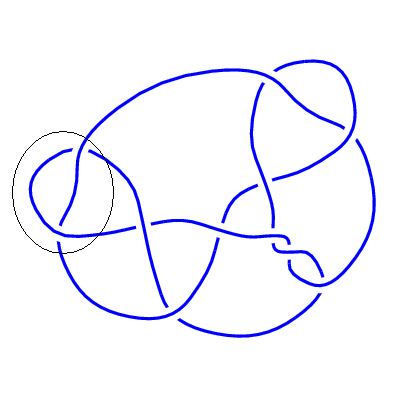}}
		\subfloat[$11n_{122}$]{\includegraphics[width=5cm]{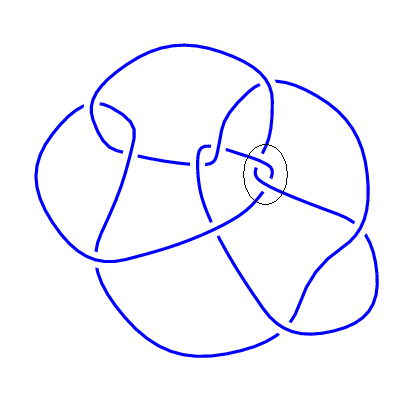}}
	\caption{Knots concordant to $3_1$}
	\label{fig:3_1}
\end{figure}

\begin{figure}[h]
		\centering
		\subfloat[$11a_{5}$]{\includegraphics[width=5cm]{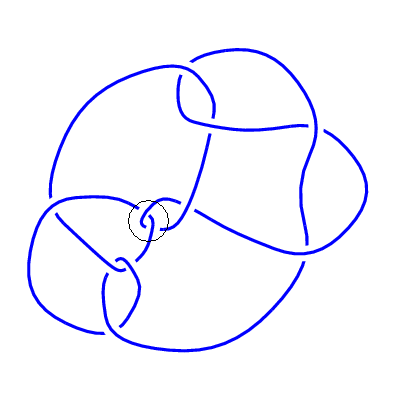}}
		\subfloat[$11a_{104}$]{\includegraphics[width=5cm]{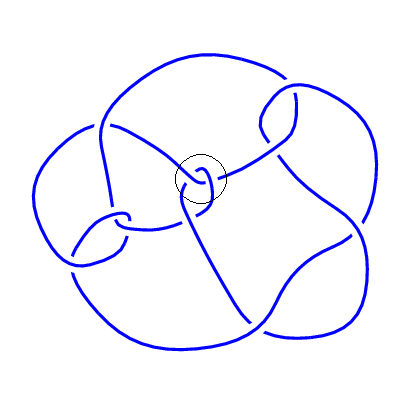}}
						
		\subfloat[$11a_{112}$]{\includegraphics[width=5cm]{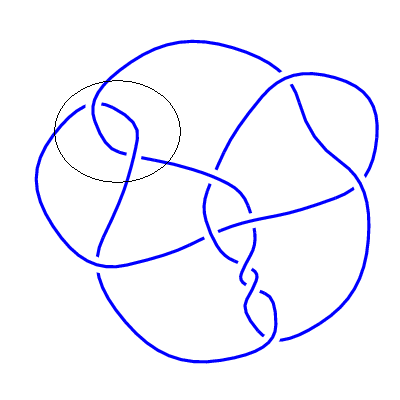}}
		\subfloat[$11a_{168}$]{\includegraphics[width=5cm]{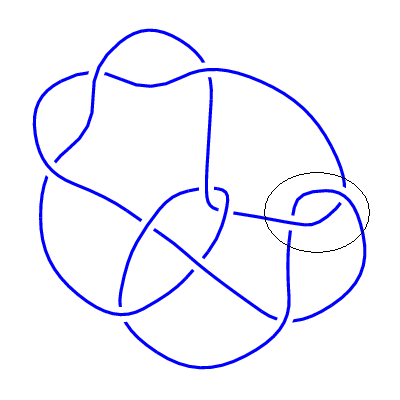}}
		
		\subfloat[$11n_{85}$]{\includegraphics[width=5cm]{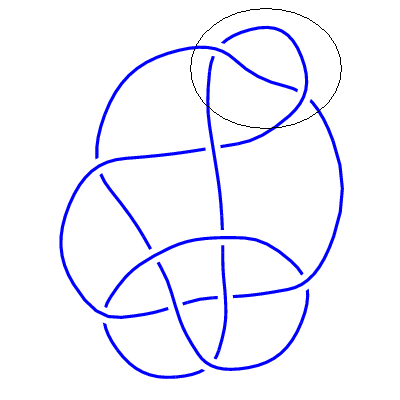}}
		\subfloat[$11n_{100}$]{\includegraphics[width=5cm]{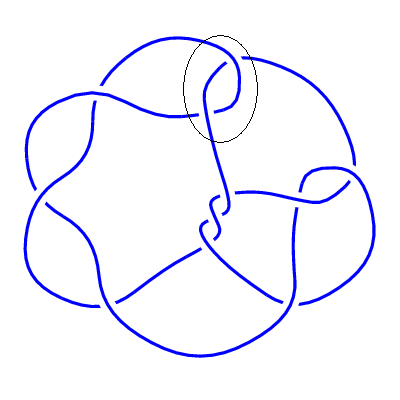}}
	\caption{Knots concordant to $4_1$}
	\label{fig:4_1}
\end{figure}

\begin{figure}[h]
		\centering
		\subfloat[$11n_{69}$]{\includegraphics[width=6cm]{11n_69.png}}
		\subfloat[$11n_{76}$]{\includegraphics[width=6cm]{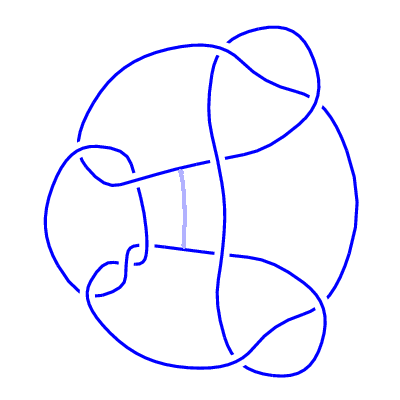}}
		
		\subfloat[$11n_{78}$]{\includegraphics[width=6cm]{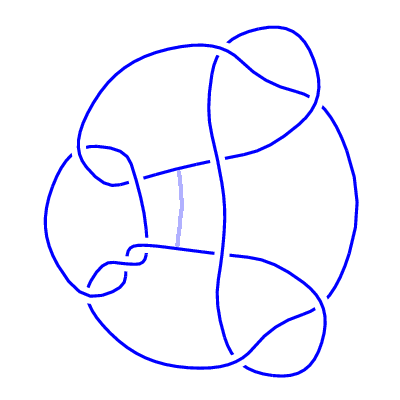}}
	\caption{Knots Concordant to $5_1$}
	\label{fig:5_1}
\end{figure}

\begin{figure}[h]
		\centering
		\subfloat[$11n_{68}$]{\includegraphics[width=6cm]{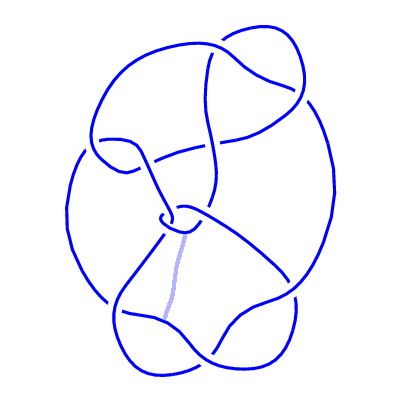}}
		\subfloat[$11n_{71}$]{\includegraphics[width=6cm]{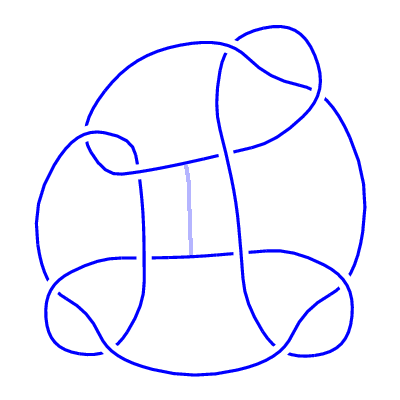}}
		
		\subfloat[$11n_{75}$]{\includegraphics[width=6cm]{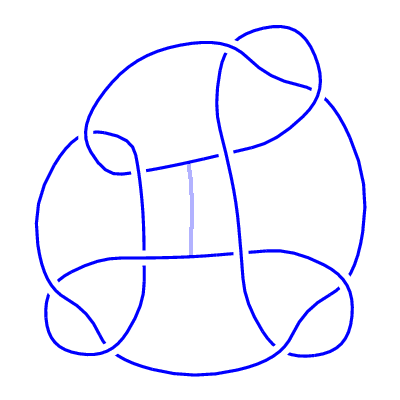}}
	\caption{Knots concordant to $5_2$}
	\label{fig:5_2}
\end{figure}

\begin{figure}[h]
		\centering
		\subfloat[$11a_{57}$]{\includegraphics[width=6cm]{11a_57.png}}
		\subfloat[$11a_{102}$]{\includegraphics[width=6cm]{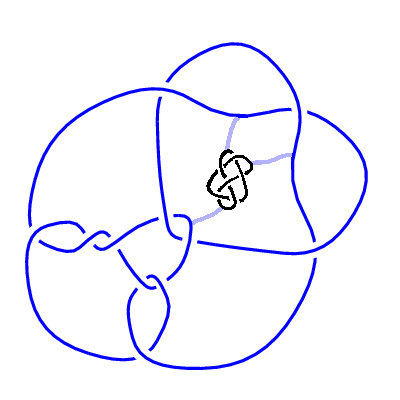}}
		
		\subfloat[$11a_{139}$]{\includegraphics[width=6cm]{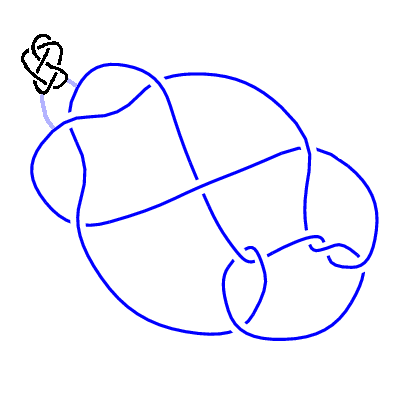}}
		\subfloat[$11a_{199}$]{\includegraphics[width=6cm]{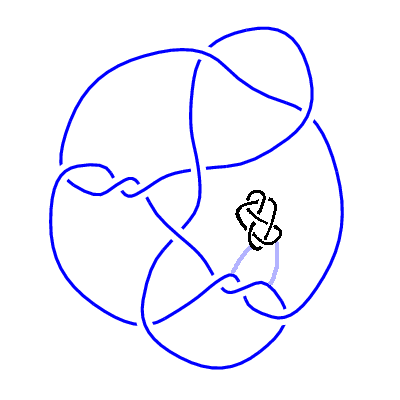}}
		
		\subfloat[$11a_{231}$]{\includegraphics[width=6cm]{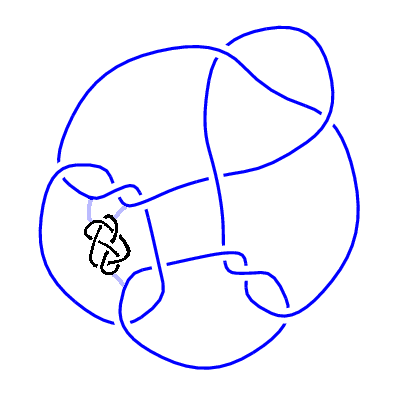}}		
	\caption{Knots concordant to $6_2$}
	\label{fig:6_2}
\end{figure}

\begin{figure}[h]
		\centering
		\subfloat[$11a_{38}$]{\includegraphics[width=6cm]{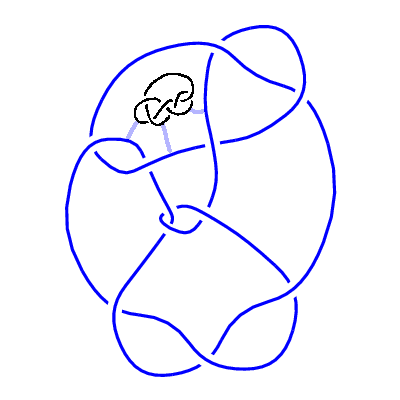}}
		\subfloat[$11a_{44}$]{\includegraphics[width=6cm]{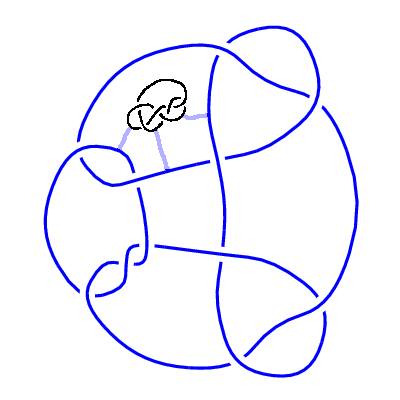}}
		
		\subfloat[$11a_{47}$]{\includegraphics[width=6cm]{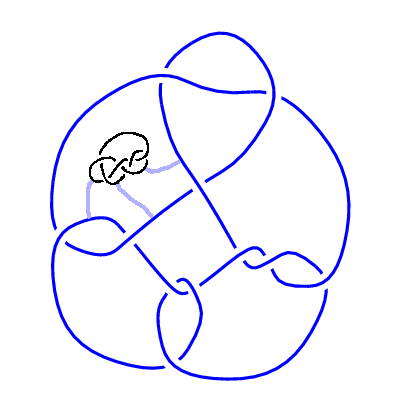}}
		\subfloat[$11a_{187}$]{\includegraphics[width=6cm]{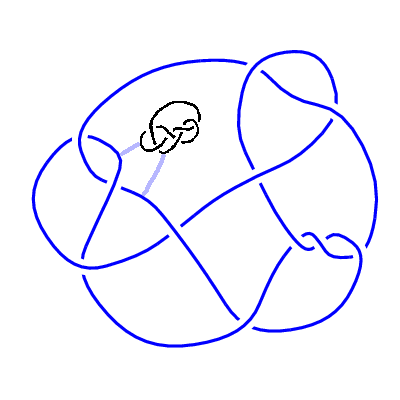}}
	\caption{Knots concordant to $6_3$}
	\label{fig:6_3}
\end{figure}

\begin{figure}[h]
		\centering
		\subfloat[$11a_{132}$]{\includegraphics[width=6cm]{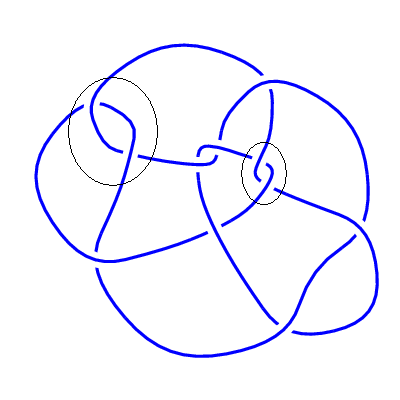}}
		\subfloat[$11a_{157}$]{\includegraphics[width=6cm]{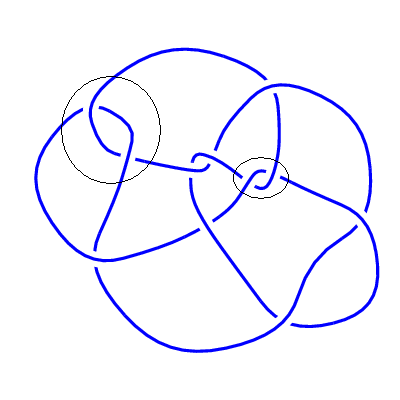}}
	\caption{Knots concordant to $3_1\#4_1$}
	\label{fig:3_1+4_1}
\end{figure}

Finally, we provide diagrams of null-concordances for the known slice 11--crossing knots.  Several of these are found to be slice by way of a concordance to $6_1$.  Since $6_1$ is slice, a concordance from a given knot to $6_1$ is enough to prove that the knot is slice.

\begin{figure}[h]
		\centering
		\subfloat[$11a_{28}$]{\includegraphics[width=.25\textwidth]{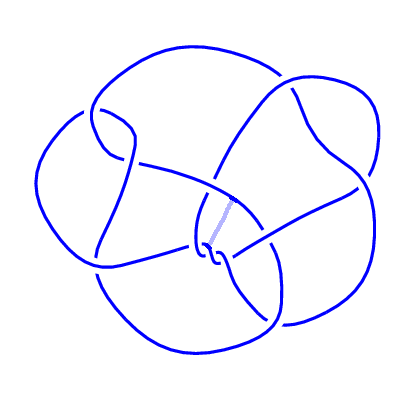}}
		\subfloat[$11a_{35}$]{\includegraphics[width=.25\textwidth]{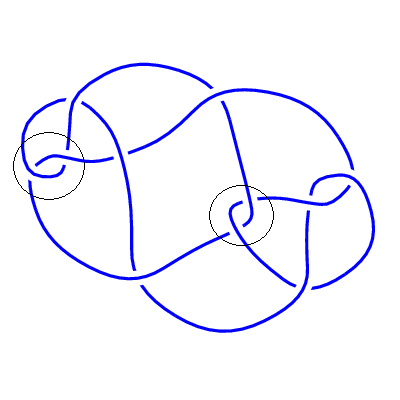}}
		\subfloat[$11a_{36}$]{\includegraphics[width=.25\textwidth]{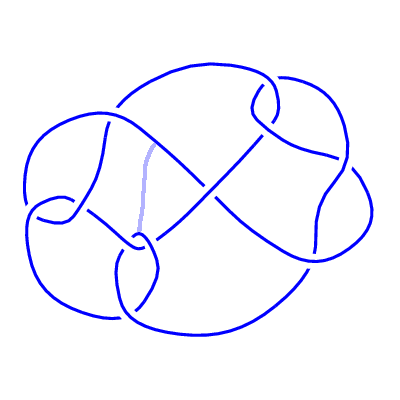}}
		
		\subfloat[$11a_{58}$]{\includegraphics[width=.25\textwidth]{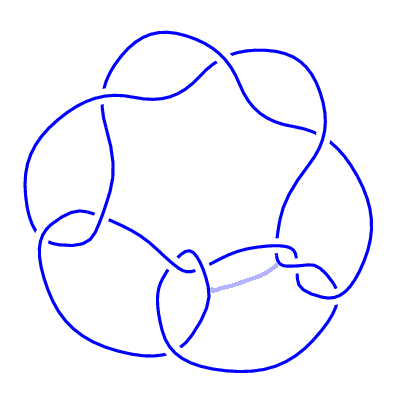}}
		\subfloat[$11a_{87}$]{\includegraphics[width=.25\textwidth]{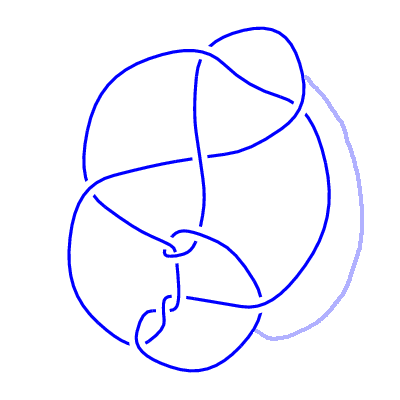}}
		\subfloat[$11a_{96}$]{\includegraphics[width=.25\textwidth]{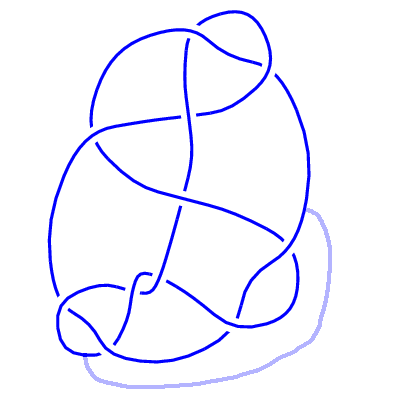}}
				
		\subfloat[$11a_{103}$]{\includegraphics[width=.25\textwidth]{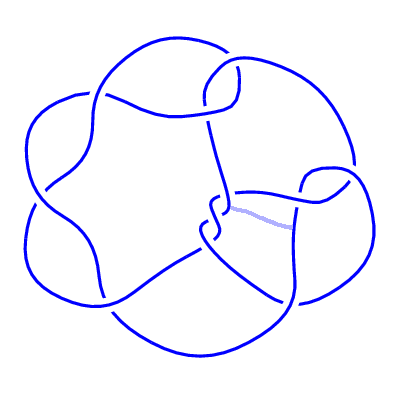}}
		\subfloat[$11a_{115}$]{\includegraphics[width=.25\textwidth]{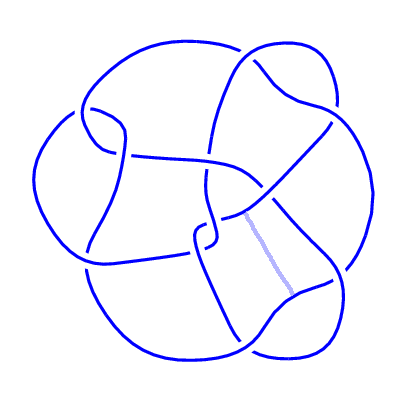}}
		\subfloat[$11a_{164}$]{\includegraphics[width=.25\textwidth]{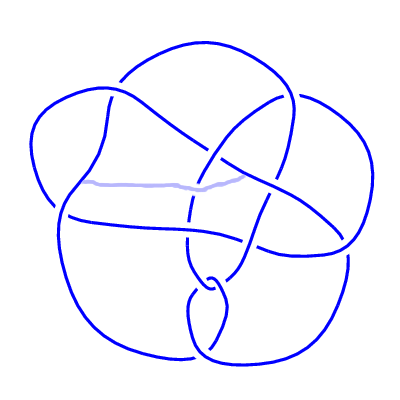}}	
		\subfloat[$11a_{165}$]{\includegraphics[width=.25\textwidth]{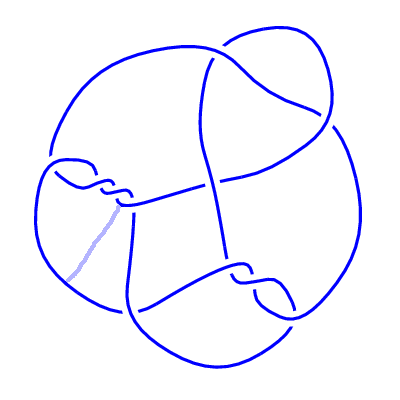}}
		
		\subfloat[$11a_{169}$]{\includegraphics[width=.25\textwidth]{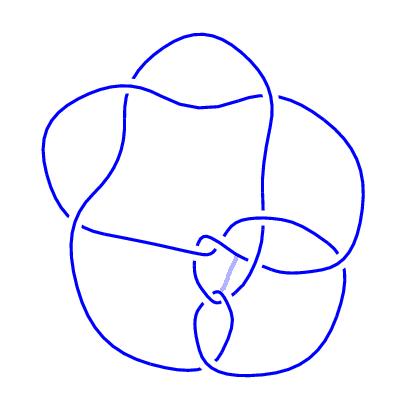}}
		\subfloat[$11a_{201}$]{\includegraphics[width=.25\textwidth]{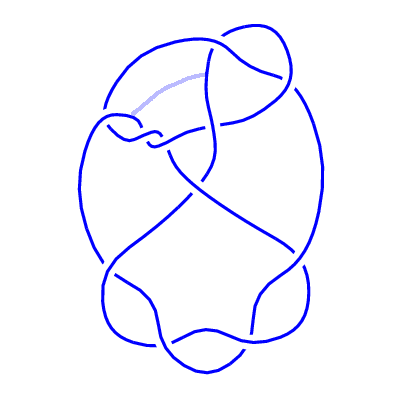}}	
		\subfloat[$11a_{316}$]{\includegraphics[width=.25\textwidth]{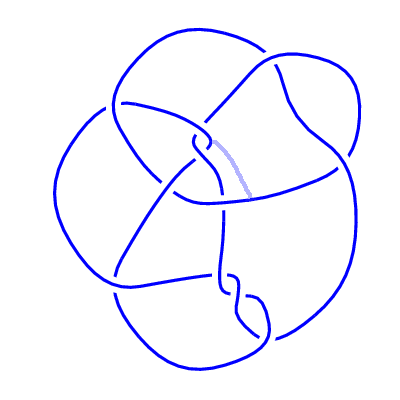}}
		\subfloat[$11a_{326}$]{\includegraphics[width=.25\textwidth]{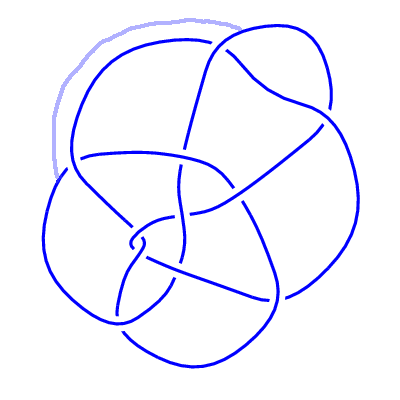}}
			\caption{Slice knots, alternating}
	\label{fig:slice}
\end{figure}

\begin{figure}[h]
		\centering
		\subfloat[$11n_{4}$]{\includegraphics[width=.25\textwidth]{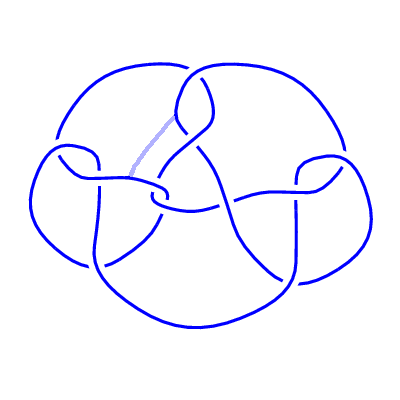}}		
		\subfloat[$11n_{21}$]{\includegraphics[width=.25\textwidth]{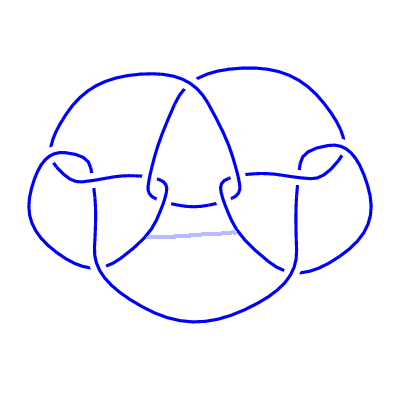}}
		\subfloat[$11n_{37}$]{\includegraphics[width=.25\textwidth]{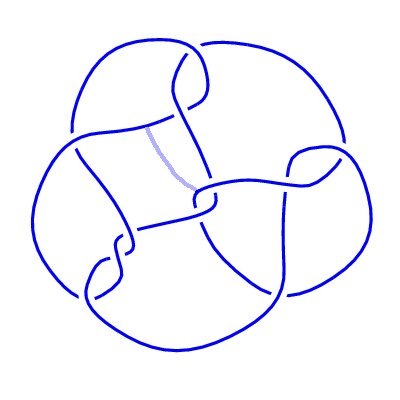}}
		\subfloat[$11n_{39}$]{\includegraphics[width=.25\textwidth]{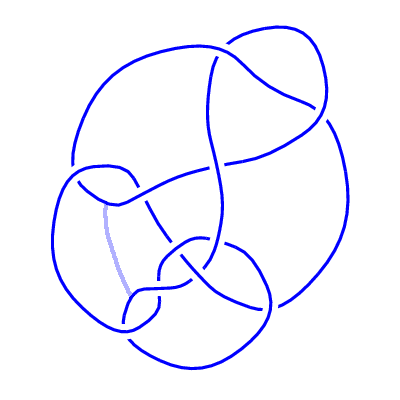}}		
		
		\subfloat[$11n_{42}$]{\includegraphics[width=.25\textwidth]{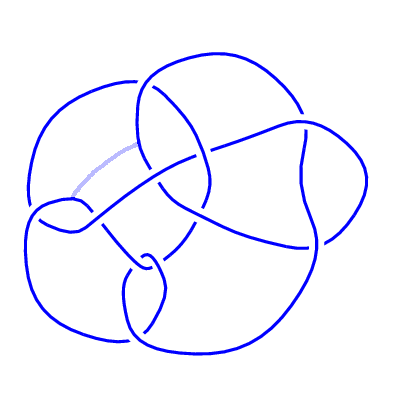}}
		\subfloat[$11n_{49}$]{\includegraphics[width=.25\textwidth]{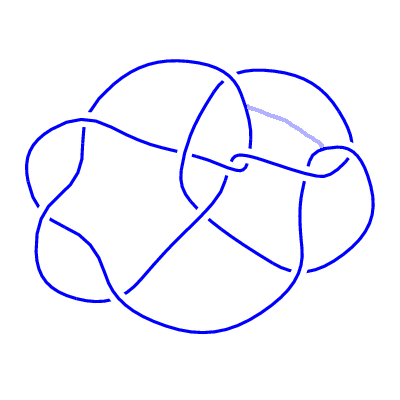}}
		\subfloat[$11n_{50}$]{\includegraphics[width=.25\textwidth]{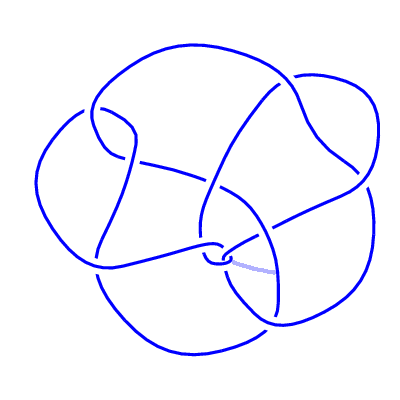}}		
		\subfloat[$11n_{67}$]{\includegraphics[width=.25\textwidth]{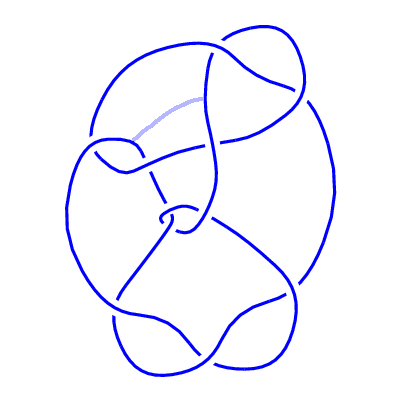}}
		
		\subfloat[$11n_{73}$]{\includegraphics[width=.25\textwidth]{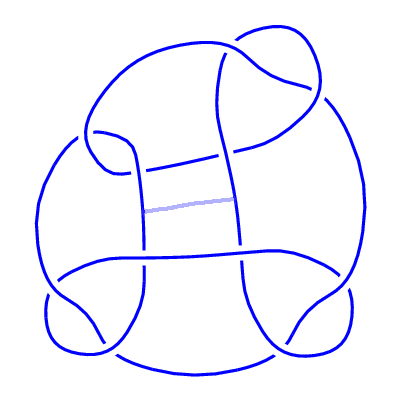}}
		\subfloat[$11n_{74}$]{\includegraphics[width=.25\textwidth]{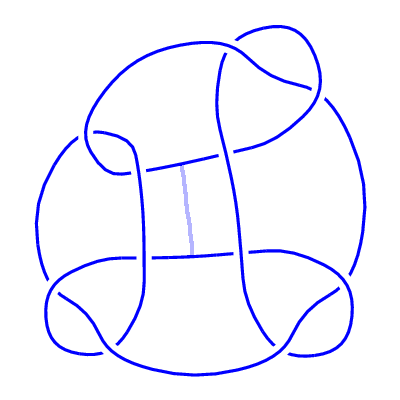}}		
		\subfloat[$11n_{83}$]{\includegraphics[width=.25\textwidth]{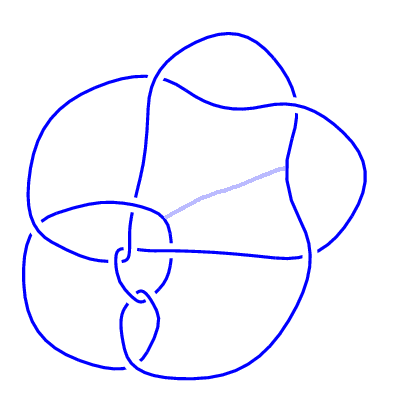}}
		\subfloat[$11n_{97}$]{\includegraphics[width=.25\textwidth]{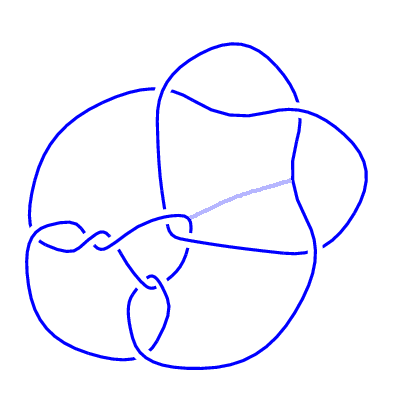}}
		
		\subfloat[$11n_{116}$]{\includegraphics[width=.25\textwidth]{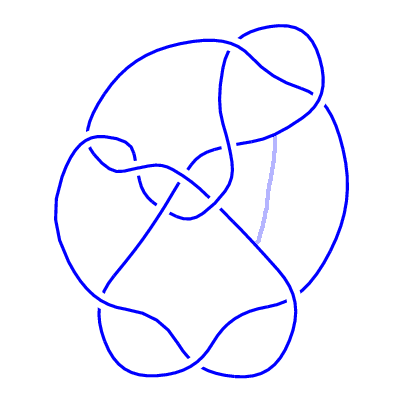}}	
		\subfloat[$11n_{132}$]{\includegraphics[width=.25\textwidth]{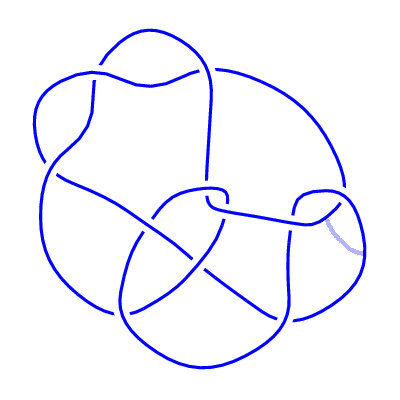}}
		\subfloat[$11n_{139}$]{\includegraphics[width=.25\textwidth]{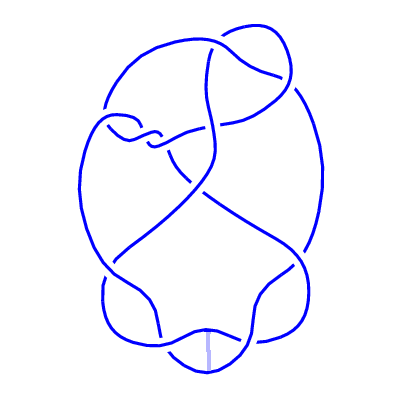}}
		\subfloat[$11n_{172}$]{\includegraphics[width=.25\textwidth]{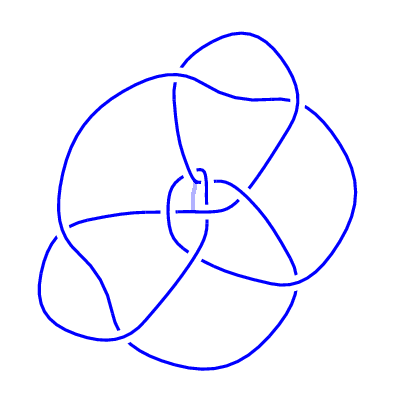}}			

	\caption{Slice knots, nonalternating}
	\label{fig:slice}
\end{figure}

\cleardoublepage

\bibliographystyle{acm}
\bibliography{MKKrefs}

 \end{document}